\documentclass{article}  

\usepackage{amssymb}
\usepackage{amsmath}
\usepackage{amsfonts}
\usepackage{amscd}
\usepackage{amsthm}

\newtheorem{definition}{Definition}

\newtheorem{example}{Example}

\newcommand{\Q}{\mathbb{Q}}
\newcommand{\Z}{\mathbb{Z}}
\newcommand{\C}{\mathbb{C}}

\newcommand{\SV}{\small \mbox{1728\,}}

\newcommand{\pFq}[5]{{}_{#1}\mathrm{F}_{#2} \left( \genfrac{}{}{0pt}{}{#3}{#4} \bigg| {#5} \right)}

\title{Computing Puiseux Expansions at Cusps of the Modular Curve $X_0(N)$}

\author{Mark van Hoeij\thanks{Supported by NSF grant 1017880.}  \\
Florida State University, Tallahassee, FL 32306-3027, USA \\
hoeij@math.fsu.edu}

\begin{document}
\maketitle
\section{Notation}


The function field of the modular curve $X_0(N)$
can be written as 
$\C(x_1)[x_2] / (\phi_N)$
where $\phi_N \in \Z[x_1,x_2]$ satisfies $\phi_N(j(\tau), j(N \tau)) = 0$.
A place $P$ corresponds to a discrete valuation $v_P$ on the function field
($v_P(g)<0$ means that a function $g$ has a pole at $P$, and $v_P(g)>0$ means a root, of that order).
The cusps are the places $P$ where $x_1,x_2$ have poles. \\ 

\noindent {\bf Goal:}  An efficient algorithm to compute Puiseux expansions at cusps of $X_0(N)$. \\

%

\noindent A Puiseux expansion at a cusp $P$ of $X_0(N)$ can be written as
\begin{equation} \label{puis}
x_1 = t^{-d}, \ \ \ \ \ x_2 = c_0 \cdot t^{-n} \cdot (1 + \cdots) \ \in \ \Z[c_0,d^{-1}][[\,t\,]].
\end{equation}
Here $n,d$ are positive integers, $c_0$ is a root of unity, $t$ is a local parameter at $P$, and
the dots refer to terms with positive powers of $t$.
To avoid negative exponents, we switch to the variables
\[ x = \frac{1}{x_1} = \frac{1}{j(\tau)}, {\rm \ \  and \ \ } h = \frac{1}{x_2} = \frac{1}{j(N \tau)}.  \]
Now $x,h$ satisfy an algebraic relation $P_N(x, h) = 0$ that is trivially obtained from $\phi_N$
by substituting $(x_1,x_2) \mapsto (x^{-1}, h^{-1})$.
However, $\phi_N$ and $P_N$ are not needed for computing a Puiseux expansion at a cusp.


In terms of $x,h$ the Puiseux expansion~(\ref{puis}) looks like
\[ x = t^d, \ \ \ h = c \cdot t^n \cdot (1+ \cdots) \]
where $c = 1/c_0$.  We can rewrite that to (from now on we will use this form):
\begin{equation} \label{puisx} h = c \cdot x^q \cdot (1 + \cdots) \ \in \ \Z[c,d^{-1}][[\,x^{1/d}\,]]
\end{equation}
where $q = \frac{n}{d}$ is a positive rational number.
We will call
\[ T := c \cdot x^q \]
the {\em initial term} of $h$.
We are only interested in those $h$ for which $\exists_N P_N(x,h)=0$.
Such $h$ turn out (see Section~\ref{SectionAlgorithm}) to be uniquely
determined by their initial term. Section~\ref{SectionN} will 
explain how to find $N$ from $T$.

\subsection{Puiseux series} \label{PuiseuxNotation}
Let
\[ \hat{K} := \bigcup_{d=1}^{\infty} \C((x^{1/d})) \]
denote the field of Puiseux series over $\C$.
If $\alpha \in \hat{K} - \{0\}$, then $v(\alpha) \in \Q$ denotes
the exponent of the initial term.  So
\[ v(h) = v(T) = q. \]
\begin{definition} \label{defprecision}
Let $a \in \frac1d \Z$ and $h$ as in~(\ref{puisx}).  By computing $h$ to {\em precision} $a$ we mean
computing the factor $(1+\cdots)$ mod $x^a$,  and hence $h$ mod $x^{q+a}$.
\end{definition}

\noindent {\bf Input and output of our algorithm.}
Given $h$ up to precision $a$, we will show that $h$ can be computed quickly to precision $2a$.
Starting with the initial value $cx^q(1 + O(x^{1/d}))$,
we will thus find $h$ mod $x^{q+d^{-1}{2^k}}$ after $k$ steps.

\subsection{The number $N$}
\label{SectionN}
For any monomial $T = c x^q$, with $c$ a root of unity, and $q$ a positive rational number, our algorithm will compute
a specified number of terms of a Puiseux series $h = T \cdot (1+\cdots)$
for which $P_N(x,h)=0$ for one $N$.
We can quickly determine $N$ from $T$.  For instance, if $N$ is prime, then either $(c,q) = (1,N)$ or $(c,q) = (\zeta_N^s, 1/N)$ for
some $s \in \{0,\ldots,N-1\}$.
The relation between other $N$'s and their $T$'s comes from composition, as shown in these examples:
\begin{example} \label{Ex1} Composing $T = \pm x^{1/2}$ and $T=x^2$  (all belonging to $N=2$) we obtain $x$ and $-x$.
Now $T = x$ belongs to $N = 1 = 2/2$, but $T = -x$ does not.
So it must belong to $N = 2 \cdot 2 = 4$.
Similarly, $i x^2$  (where $i = \zeta_4 = \sqrt{-1}$) belongs to $X_0(2^k)$ for some $k$   
since it can be obtained by repeated compositions of $x^2$ and $\pm x^{1/2}$.  Here $k$ must be 5 since we do not obtain
$i x^2$ by composing fewer than 5 functions from $\{x^2, \pm x^{1/2}\}$. In contrast, $T = i x^{1/2}$ belongs to $N=2^3$. \\
Likewise, $-x^{1/3}$ and $x^{3/4}$ can only belong to $N=3 \cdot 2^2 = 12$, because both require compositions involving 1 element
from $\{x^3, \zeta_3^*  x^{1/3}\}$ and 2 elements from $\{x^2, \pm x^{1/2}\}$.
Likewise,  $T=x^{5/3}$ and $T = \zeta_3 x^{5/3}$ belong to $N=15$,  and $T = \zeta_5 x$ belongs to $N = 25$.
\end{example}




\section{A relation between $x$ and $h$}
The reciprocals of $x$ and $h$ satisfy the modular equation $\phi_N$.
Since $\phi_N$ can be large when $N$ is large, we will use another relation between $x$ and $h$,  one that is valid for any $N$.
Define $E, F, G \in \Z[[x]]$ as
\[ E := x \cdot \sqrt{1-\SV x}, \ \ \ F := \pFq{2}{1}{\frac{1}{12},\frac{5}{12}}{1}{\SV x} \]
and
\[ G := E \cdot F^2 =  E \cdot \pFq{3}{2}{\frac16, \frac12, \frac56}{1,1}{\SV x}. \]
$G$ satisfies a linear homogeneous differential equation $L_3$ over $\Q(x)$
\[ L_3: \ \ G''' + a_2 G'' + a_1 G' + a_0 G = 0 \]
with $a_0,a_1,a_2 \in \Q(x)$.
The factor $E$ in $G$ was selected to ensure $a_2 = 0$.

The following relation
\begin{equation} \label{main} v(h) \cdot G \circ h = h' \cdot G \end{equation}
holds for every $h \in \hat{K}$ for which $\exists_N  P_N(x,h) = 0$.
We computed this relation by reformulating the condition [To Do: Find reference]
that $F \circ h$ should be an algebraic function times $F$.
The projective monodromy matrices of $F$ are precisely the famous generators of the modular group PSL(2,$\Z$).

Section~\ref{SectionAlgorithm} shows that for any positive rational number $q$ and any $c \in \C - \{0\}$
there exists precisely one $h \in \hat{K}$ that satisfies~(\ref{main}) and has $c x^q$ as its initial term.
This $h$ is algebraic over $\C(x)$ iff\footnote{The fact that $h$ satisfies some $P_N$ when $c$ is a root of unity
implies that $h$ can not be algebraic when $c$ is not a root of unity.
If $c$ is not a root of unity, and if $h$ were algebraic,
then $c \equiv \zeta_N$ mod $p$ for a large $N$ and a large prime $p$, and we would
get arbitrarily high lower bounds on the algebraic degree of $h$ reduced mod $p$, leading to a contradiction.}
$c$ is a root of unity.



\section{Computing $h$ from its initial term}
\label{SectionAlgorithm}
Differentiating~(\ref{main}) and dividing by $h'$ we find
\begin{equation} \label{formule2}
v(h) \cdot  G' \circ h  = G \cdot (\frac{h''}{h'} + \frac{G'}{G}) = G \cdot {\rm ld}(h'G) \end{equation}
where ${\rm ld}$ denotes the logarithmic derivative,  ${\rm ld}(u) := {\rm ln}(u)' = u'/u$.
Suppose that $h_0$ is an approximation of $h$ with $v(\epsilon) \geqslant a+v(h) > a$ where $\epsilon$ denotes $h - h_0$.
Substituting $h = h_0 + \epsilon$ in~(\ref{main}) and~(\ref{formule2}) gives
\begin{equation} \label{stap3}
 (h_0' + \epsilon') \cdot G = v(h) \cdot G \circ (h_0 + \epsilon) = v(h) \cdot (G \circ h_0 + \epsilon \cdot (G' \circ h_0) + O(\epsilon^2)) \end{equation}
and, using $v(G)=1$,
\begin{equation} \label{stap4}
 v(h) \cdot G' \circ h_0  \, + O(\epsilon) \  = \ G \cdot {\rm ld}(h_0' G) \, + O(x^a). \end{equation}
Substituting~(\ref{stap4}) into~(\ref{stap3}), dividing by $G$, then subtracting $h_0'$, gives
\[ \epsilon' = \frac{v(h) \cdot (G \circ h_0)}{G} + \epsilon \cdot {\rm ld}(h_0' G)  - h_0' \, + O(x^{v(h)+2a-1}). \]
Now $\epsilon' = {\rm ld}(A) \epsilon + B + O(x^{s})$ has a solution $\epsilon = A \int B/A \, + O(x^{s+1})$, applying that
gives
\begin{equation} \label{eps}
\epsilon = h_0' G \int \frac{1}{G} \left( \frac{v(h) \cdot (G \circ h_0)}{h_0'G} - 1\right) {\rm d}x \, + O(x^{v(h)+2a}).
\end{equation}
Adding this to $h_0$ doubles the precision in the sense of Definition~\ref{defprecision}. \\[5pt]

\noindent {\bf Algorithm PuiseuxX0N}. \\[5pt]
{\bf Input:} $T = c x^q$ where $c$ is a root of unity and $q$ a positive rational number, and a positive integer $k$. \\[5pt]
{\bf Output:} An approximation of precision $d^{-1} 2^k$ (as in Definition~\ref{defprecision})
of a Puiseux series $h$ with initial term $T$ that satisfies $P_N(x,h)=0$ (with $N$ as in Example~\ref{Ex1}). \\[5pt]
{\bf Step 1.}  $h_0 := T$  and  $a := d^{-1}$  where $d = {\rm denominator}(q)$. \\
{\bf Step 2.}  Repeat $k$ times: \\
\mbox{} \ \ (a) Compute $\epsilon$ mod $x^{q+2a}$ with formula~(\ref{eps}). \\
\mbox{} \ \ (b) $h_0 := h_0 + \epsilon$  and  $a := 2a$. \\
{\bf Step 3.} Return $h_0$. \\[5pt]

%
A Maple implementation is given at {\tt www.math.fsu.edu/$\sim$hoeij/files/X0N}, in the file {\tt PuiseuxX0N}.
The CPU time is dominated by the cost of composing $G \circ h_0$.  Now $G$ contains $F^2$, so we
must compose a $\mbox{}_2F_1$ function with a truncated power series $h_0$.
Brent and Kung \cite{BK} described an algorithm that can perform this step efficiently.  This,
combined with fast arithmetic in $\Z[c,d^{-1}]$, reduces the computational complexity
to quasi-linear time (logarithmic factors times the size of the output).

One could compute $\phi_N$ by (i) computing Puiseux expansions to sufficient precision, and then (ii) reconstructing $\phi_N$ from them.
Step~(i) is quasi-linear, and so is Step~(ii) if $N$ is for example a power of 2. But if $N$ contains large prime(s), it is not clear if Step~(ii)
can be done faster than \cite{POLY}.

\end{document}